\documentclass[11pt,a4paper]{article}

\usepackage[utf8]{inputenc}
\usepackage{graphicx}
\usepackage{tikz}
\usepackage{booktabs}
\usepackage{subcaption}
\usepackage{pdflscape}
\usepackage{array}
\usetikzlibrary{arrows.meta,positioning}

\usepackage{amsmath, amsthm, amssymb}
\usepackage{algorithm}
\usepackage[noend]{algpseudocode}
\usepackage{microtype}
\usepackage{hyperref}
\newtheorem{remark}{Remark}

\newtheorem{property}{Property} 
 
\newtheorem{lemma}{Lemma} 

\begin{document}

\title{Revisiting Johnson's rule for minimizing makespan in the two-machine flow shop scheduling problem}

\author{Federico Della Croce$^{1,2}$, Quentin Schau$^{1,3}$ \\
	$^1$DIGEP Politecnico di Torino, Italy,
	\\
	$^2$ CNR, IEIIT, Torino, Italy \\
	\texttt{federico.dellacroce@polito.it}\\
	$^3$LIFAT University of Tours, France,\\
	\texttt{quentin.schau@polito.it}}
\date{}
\maketitle

%{\footnotesize
\begin{abstract}
We consider Johnson’s rule for minimizing the makespan in the two-machine flow shop problem.
Although its worst-case time complexity is O(n logn), 
we show that it is possible to detect in linear time whether a full sorting of jobs can be avoided and an optimal solution can be computed in O(n) time.
A probabilistic analysis indicates that linear-time complexity holds with high probability under uniformly distributed processing times, 
a result further supported by extensive computational experimentation.
\end{abstract}

\medskip

\noindent
{\bf Keywords}: two-machine flow shop scheduling; Johnson's algorithm, linear time complexity, probabilistic analysis.
%}
\section{Introduction}

We consider the classical two-machine flow shop scheduling problem with the objective of
minimizing the makespan, commonly denoted by $F2 \; || \; C_{\max}$ \cite{Graham1979}.
A set $J=\{1,\dots,n\}$ of jobs must be processed on two machines $M_1$ and $M_2$ in this
order.
Job $j$ requires a processing time $p_{1,j}$ on $M_1$ and $p_{2,j}$ on $M_2$.
All jobs are available at time~0, each machine can process at most one job at a time, and
preemption is not allowed.
It is well known that there always exists an optimal schedule that is a permutation
schedule and can be obtained by Johnson’s rule \cite{Johnson1954}.
This rule first schedules the set $A$ of jobs such that $p_{1,j} \le p_{2,j}$ in
nondecreasing order of $p_{1,j}$, and then the set $B$ of jobs such that
$p_{1,j} > p_{2,j}$ in nonincreasing order of $p_{2,j}$.

\begin{algorithm}[H]
	\caption{Johnson's rule for $F2\|\!C_{\max}$}
	\label{alg:johnson}
	\begin{algorithmic}[1]
		\State \textbf{Input:} Processing times $\{(p_{1,j},p_{2,j})\}_{j=1}^n$
		\State $A \gets \{j : p_{1,j} \le p_{2,j}\}$,\quad
		       $B \gets \{j : p_{1,j} > p_{2,j}\}$
		\State Sort $A$ by nondecreasing $p_{1,j}$; sort $B$ by nonincreasing $p_{2,j}$
		\State \Return sequence $\pi \gets$ \textsc{Concatenate}$(A,B)$
	\end{algorithmic}
\end{algorithm}

\noindent
Due to the required sorting operations, Johnson’s rule runs in $O(n \log n)$ time and
$O(n)$ space.
As noted in \cite{Pinedo2016}, schedules obtained by Johnson’s rule
``are by no means the only schedules that are optimal for $F2 \; || \; C_{\max}$.
The class of optimal schedules appears to be hard to characterize and data dependent.''
Motivated by this observation, several studies have focused on the generation of
equivalent optimal solutions.
The characterization and enumeration of all optimal sequences was first addressed
in \cite{Billaut1998}.
This work was later extended in \cite{Briand2006}, where a new sufficient condition
for optimality was proposed.
More recently, refined optimality conditions and additional structural insights into
the set of optimal schedules were introduced in \cite{Hadda2018}.

In this paper, we investigate general conditions under which the problem can be solved
optimally in linear time.
More specifically, we prove that a simple procedure running in $O(n)$ time allows one
to determine whether a full sort can be avoided and, correspondingly, whether an
optimal solution can be computed in linear time.
In this case, the procedure also identifies a large set of equivalent optimal
sequences.
To the best of our knowledge, this is the first work showing that Johnson’s rule can be
implemented in linear time with high probability under standard distributions of
processing times.

\medskip

\noindent
The remainder of the paper is organized as follows.
In Section~\ref{equiv}, we show how to detect in linear time subsets of jobs that do not
need to be sorted in the construction of an optimal solution, as any permutation of
these jobs also induces an optimal schedule.
Accordingly, we verify whether the number of remaining jobs to be sorted is small
enough to guarantee an overall linear-time complexity.
As a counterpart, we also present a family of instances for which a full sort is
necessary.
Section~\ref{prob} provides a probabilistic analysis for instances with uniformly
distributed processing times, showing that linear-time complexity is guaranteed with
probability lower bounded by $1- 10^{-3}$ for instances of size $n = 100$ and by $1-10^{-9}$ for instances of size $n = 200$.
In Section~\ref{compres}, computational experiments indicate that, for standard uniform
distributions and instance sizes $n = 100, 1000, 10000$, and $100000$, an optimal
solution is always obtained in linear time.
The same behavior is observed for instances of size $n = 1000$ generated according to
other widely used distributions.
Finally, Section~\ref{conclus} concludes the paper.

\section{Revisiting Johnson's rule for problem $F2 \; || \; C_{\max}$}
\label{equiv}

Assume, without loss of generality, that the $n = n_A + n_B$ jobs are indexed according to
Johnson’s sequence, where jobs $1,\dots,n_A$ belong to set $A$ and jobs $n_A+1,\dots,n$
belong to set $B$.
Let $p^{\max}$ denote the largest processing time among all operations.
Moreover, let $p^{\max}_{A1}$ (resp., $p^{\max}_{B2}$) denote the largest processing time
in subset $A$ on machine $M_1$ (resp., in subset $B$ on machine $M_2$).
Let
\[
P_1 = \sum_{i=1}^{n} p_{1,i}
\qquad\text{and}\qquad
P_2 = \sum_{i=1}^{n} p_{2,i}
\]
be the total processing times on machines $M_1$ and $M_2$, respectively.
Similarly, let
\[
P_{1A} = \sum_{i=1}^{n_A} p_{1,i}, \qquad
P_{2A} = \sum_{i=1}^{n_A} p_{2,i}
\]
be the total processing times of the jobs in subset $A$ on machines $M_1$ and $M_2$, and let
\[
P_{1B} = \sum_{i=n_A+1}^{n} p_{1,i}, \qquad
P_{2B} = \sum_{i=n_A+1}^{n} p_{2,i}
\]
be the corresponding totals for subset $B$.
Finally, for each job $i$, let $C_{1,i}$ (resp., $C_{2,i}$) denote its completion time on
machine $M_1$ (resp., on machine $M_2$).

\medskip

\begin{remark}
The problem $F2 \; || \; C_{\max}$ satisfies the reversibility property stated below.
Consider a pair of two-machine flow shops where $p^{(1)}_{i,j}$ and $p^{(2)}_{i,j}$ denote the
processing time of job $j$ on machine $i$ in the first and second flow shop, respectively.
\begin{property}[\cite{Muth1979,Pinedo2016}]
\label{prop11}
If $p^{(1)}_{1,j}=p^{(2)}_{2,j}$ and $p^{(1)}_{2,j}=p^{(2)}_{1,j}$ for all jobs $j$, then
sequence $j_1,\dots,j_n$ in the first flow shop results in the same makespan as
sequence $j_n,\dots,j_1$ in the second flow shop.
\end{property}
\end{remark}

\noindent
The following properties provide sufficient conditions for finding an optimal solution
to $F2 \; || \; C_{\max}$ without requiring a full sort of the job set.

\begin{property}
\label{prop1}
If \( P_{1A} \leq P_{2A} - p^{\max}_{A1} \), let \( k_A \) denote the smallest index such that\\
\[\sum_{i \in A: \, p_{1,i}\leq p_{1,k_A}} p_{1,i} \leq
\sum_{i \in A: \, p_{1,i}\leq p_{1,k_A}} p_{2,i} - p^{\max}_{A1}\]
and
\[p_{1,k_A} < p_{1,k_A+1}.\]
Then, for set \( A \), an optimal schedule is obtained by sequencing jobs
\(1,\dots,k_A\) as in Johnson’s order, followed by jobs \(k_A+1,\dots,n_A\) in any order.
\end{property}

\begin{proof}
Since \(p_{1,i}\leq p_{2,i}\) for all \(i\in A\), once
\[
\sum_{j=1}^{k_A} p_{1,j} \leq \sum_{j=1}^{k_A} p_{2,j} - p^{\max}_{A1},
\]
then for any \(i^*\in\{k_A+1,\dots,n_A\}\) placed in position \(k_A+1\) we have
\[
C_{1,i^*}=\sum_{j=1}^{k_A} p_{1,j} + p_{1,i^*}
\leq
\sum_{j=1}^{k_A} p_{1,j} + p^{\max}_{A1}
\leq
\sum_{j=1}^{k_A} p_{2,j},
\]
and hence no idle time is incurred on $M_2$ for job \(i^*\).
Moreover, from position \(k_A+1\) onward, the gap between completion times on the two
machines cannot decrease with respect to the gap at position \(k_A\); therefore, no idle
time can occur on $M_2$ for positions \(k_A+1,\dots,n_A\).
Consequently, the completion time of the last job of set \(A\) on $M_2$ is the same for any
permutation of the jobs in positions \(k_A+1,\dots,n_A\).
The additional condition \(p_{1,k_A} < p_{1,k_A+1}\) ensures that all jobs in $A$ with
$p_{1,i}=p_{1,k_A}$ are included in the portion that is not declared freely permutable.
\end{proof}

\begin{property}
\label{prop2}
If \( P_{2B} \leq P_{1B} - p^{\max}_{B2} \), let \( k_B \) denote the largest index such that
\[
\sum_{i \in B: \, p_{2,i}\leq p_{2,k_B}} p_{2,i}
\leq
\sum_{i \in B: \, p_{2,i}\leq p_{2,k_B}} p_{1,i} - p^{\max}_{B2}
\]
and
\[
p_{2,k_B} > p_{2,k_B-1}.
\]
Then, for set \(B\), an optimal schedule is obtained by sequencing jobs
\(n_A+1,\dots,k_B-1\) in any order, followed by jobs \(k_B,\dots,n\) as in Johnson’s order.
\end{property}

\begin{proof}
By the reversibility property (Property~\ref{prop11}), the conditions of
Property~\ref{prop2} correspond to the conditions of Property~\ref{prop1} applied to the
reversed instance.
\end{proof}

\noindent
Combining Properties~\ref{prop1} and~\ref{prop2}, we can decide in linear time whether a
sufficient condition for solving $F2 \; || \; C_{\max}$ in overall linear time holds.

\begin{property}
\label{coro1}
If Properties~\ref{prop1} and~\ref{prop2} hold and
$k_A \log k_A \leq n$ and $(n-k_B+1)\log(n-k_B+1)\leq n$, then an optimal solution to
$F2 \; || \; C_{\max}$ can be computed in linear time:
only jobs $1,\dots,k_A$ need to be sorted within subset $A$, and only jobs $k_B,\dots,n$
need to be sorted within subset $B$.
\end{property}

\begin{proof}
Partitioning the jobs into subsets $A$ and $B$ takes linear time.
Finding $k_A$ in $A$ and $k_B$ in $B$ can be done in linear time via selection algorithms,
e.g., the median-finding technique of \cite{Blum1972} (similarly to the identification of
the critical item in the continuous relaxation of the $0/1$ knapsack problem
\cite{BZ80,MT90}).
If $k_A \log k_A \leq n$ and $(n-k_B+1)\log(n-k_B+1)\leq n$, then sorting the reduced sets
$1,\dots,k_A$ and $k_B,\dots,n$ also requires overall linear time.
\end{proof}

\begin{figure*}
	\centering

	\renewcommand{\arraystretch}{1.2}
	\setlength{\tabcolsep}{5pt}
	{\footnotesize
		\begin{tabular}{@{}l*{14}{c}@{}}
			\toprule
			Job $j$   & 1 & 2 & 3 & 4 & 5 & 6 & 7 & 8 & 9 & 10 & 11 & 12 & 13 & 14 \\
			\midrule
			$p_{1,j}$ & 1 & 2 & 3 & 4 & 5 & 6 & 7 & 9 & 8 & 8 & 8 & 7 & 9 & 10 \\
			$p_{2,j}$ & 8 & 9 & 7 & 8 & 9 & 7 & 9 & 7 & 6 & 5 & 4 & 4 & 3 & 2 \\
			\bottomrule
		\end{tabular}
	}

		\hspace*{-1cm}
	\begin{tikzpicture}[scale=0.15, every node/.style={font=\tiny}]
		\draw[->] (0,0) -- (92,0) node[right]{Time};
		\draw (0,6) node[left]{$M_1$};
		\draw (0,3) node[left]{$M_2$};

		\begin{scope}[xscale=0.9]
			\foreach \s/\d/\id in {
				0/1/1,
				1.111/2/2,
				3.333/3/3,
				6.667/4/4,
				11.111/5/5,
				16.667/6/6,
				23.333/7/7,
				31.111/9/8,
				41.111/8/9,
				50/8/10,
				58.889/8/11,
				67.778/7/12,
				75.556/9/13,
				85.556/10/14
			}{
				\ifnum\id<3
					\draw[fill=red!40] (\s,5) rectangle ++(\d,2);
				\else\ifnum\id<8
					\draw[fill=green!40] (\s,5) rectangle ++(\d,2);
				\else\ifnum\id<13
					\draw[fill=blue!40] (\s,5) rectangle ++(\d,2);
				\else
					\draw[fill=red!40] (\s,5) rectangle ++(\d,2);
				\fi\fi\fi
			}

			\foreach \x in {0,1.111,3.333,6.667,11.111,16.667,23.333,31.111,41.111,50,58.889,67.778,75.556,85.556,96.667}{
				\draw[thin] (\x,5) -- (\x,7);
			}

			\foreach \s/\d/\id in {
				1.111/8/1,
				10/9/2,
				20/7/3,
				27.778/8/4,
				36.667/9/5,
				46.667/7/6,
				54.444/9/7,
				64.444/7/8,
				72.222/6/9,
				78.889/5/10,
				84.444/4/11,
				88.889/4/12,
				93.333/3/13,
				96.667/2/14
			}{
				\ifnum\id<3
					\draw[fill=red!40] (\s,2) rectangle ++(\d,2);
				\else\ifnum\id<8
					\draw[fill=green!40] (\s,2) rectangle ++(\d,2);
				\else\ifnum\id<13
					\draw[fill=blue!40] (\s,2) rectangle ++(\d,2);
				\else
					\draw[fill=red!40] (\s,2) rectangle ++(\d,2);
				\fi\fi\fi
			}

			\draw[dashed,red] (98.889,0) -- (98.889,7) node[above]{\small \hspace*{0.7cm} $C_{\max}=89$};
		\end{scope}

		\foreach \s/\d/\id in {
			0/1/1,
			1/2/2,
			3/3/3,
			6/4/4,
			10/5/5,
			15/6/6,
			21/7/7,
			28/9/8,
			37/8/9,
			45/8/10,
			53/8/11,
			61/7/12,
			68/9/13,
			77/10/14
		}{
			\node[font=\tiny] at ({\s+\d/2},7.9) {\id};
		}

		% Sotto M2
		\foreach \s/\d/\id in {
			1/8/1,
			9/9/2,
			18/7/3,
			25/8/4,
			33/9/5,
			42/7/6,
			49/9/7,
			58/7/8,
			65/6/9,
			71/5/10,
			76/4/11,
			80/4/12,
			84/3/13,
			87/2/14
		}{
			\node[font=\tiny] at ({\s+\d/2},1.1) {\id};
		}

		\node[below right] at (0,0) {\small Sequence: $1,2,3,4,5,6,7,8,9,10,11,12,13,14$};
	\end{tikzpicture}

	\bigskip

	\hspace*{-1cm}
	\begin{tikzpicture}[scale=0.15, every node/.style={font=\tiny}]
		\draw[->] (0,0) -- (92,0) node[right]{Time};
		\draw (0,6) node[left]{$M_1$};
		\draw (0,3) node[left]{$M_2$};

		\begin{scope}[xscale=0.9]
			\foreach \s/\d/\id in {
				0/1/1,
				1.111/2/2,
				3.333/7/7,
				11.111/6/6,
				17.778/5/5,
				23.333/4/4,
				27.778/3/3,
				31.111/7/12,
				38.889/8/11,
				47.778/8/10,
				56.667/8/9,
				65.556/9/8,
				75.556/9/13,
				85.556/10/14
			}{
				\ifnum\id<3
					\draw[fill=red!40] (\s,5) rectangle ++(\d,2);
				\else\ifnum\id<8
					\draw[fill=green!40] (\s,5) rectangle ++(\d,2);
				\else\ifnum\id<13
					\draw[fill=blue!40] (\s,5) rectangle ++(\d,2);
				\else
					\draw[fill=red!40] (\s,5) rectangle ++(\d,2);
				\fi\fi\fi
			}

			\foreach \x in {0,1.111,3.333,11.111,17.778,23.333,27.778,31.111,38.889,47.778,56.667,65.556,75.556,85.556,96.667}{
				\draw[thin] (\x,5) -- (\x,7);
			}

			\foreach \s/\d/\id in {
				1.111/8/1,
				10/9/2,
				20/9/7,
				30/7/6,
				37.778/9/5,
				47.778/8/4,
				56.667/7/3,
				64.444/4/12,
				68.889/4/11,
				73.333/5/10,
				78.889/6/9,
				85.556/7/8,
				93.333/3/13,
				96.667/2/14
			}{
				\ifnum\id<3
					\draw[fill=red!40] (\s,2) rectangle ++(\d,2);
				\else\ifnum\id<8
					\draw[fill=green!40] (\s,2) rectangle ++(\d,2);
				\else\ifnum\id<13
					\draw[fill=blue!40] (\s,2) rectangle ++(\d,2);
				\else
					\draw[fill=red!40] (\s,2) rectangle ++(\d,2);
				\fi\fi\fi
			}

			\draw[dashed,red] (98.889,0) -- (98.889,7) node[above]{\small \hspace*{0.7cm} $C_{\max}=89$};
		\end{scope}

		\foreach \s/\d/\id in {
			0/1/1,
			1/2/2,
			3/7/7,
			10/6/6,
			16/5/5,
			21/4/4,
			25/3/3,
			28/7/12,
			35/8/11,
			43/8/10,
			51/8/9,
			59/9/8,
			68/9/13,
			77/10/14
		}{
			\node[font=\tiny] at ({\s+\d/2},7.9) {\id};
		}

		% Sotto M2
		\foreach \s/\d/\id in {
			1/8/1,
			9/9/2,
			18/9/7,
			27/7/6,
			34/9/5,
			43/8/4,
			51/7/3,
			58/4/12,
			62/4/11,
			66/5/10,
			71/6/9,
			77/7/8,
			84/3/13,
			87/2/14
		}{
			\node[font=\tiny] at ({\s+\d/2},1.1) {\id};
		}

		\node[below right] at (0,0) {\small Sequence: $1,2,7,6,5,4,3,12,11,10,9,8,13,14$};
	\end{tikzpicture}

	\caption{\small 14-job instance: Gantt diagram of Johnson's schedule and of an equivalent optimal schedule.}
	\label{fff1}
\end{figure*}
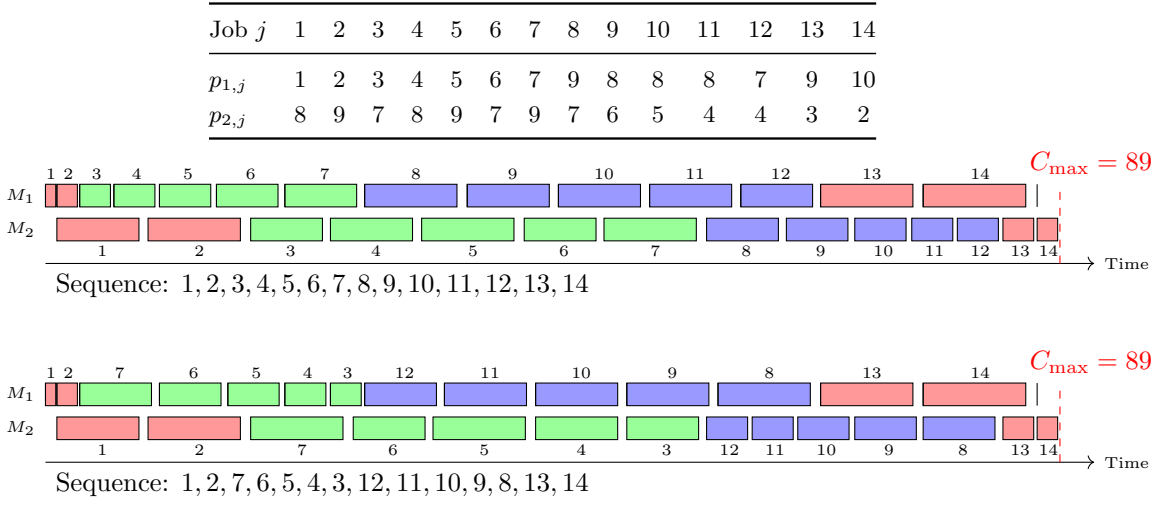

\noindent 
Consider the instance in Figure~\ref{fff1} with \(14\) jobs, for which
Properties~\ref{prop1} and~\ref{prop2} hold with \(k_A = n-k_B = 2\).
Accordingly, in an optimal solution, for subset \(A\) only the first two jobs (shown in red)
need to be sorted and placed at the beginning of the sequence.
Similarly, for subset \(B\) only the last two jobs (also shown in red) need to be sorted and
placed at the end of the sequence.
\medskip

\noindent
For the remaining jobs \(\{3,4,5,6,7\}\) in subset \(A\) (shown in green), any 
permutation yields an optimal solution.
The same holds for the remaining jobs \(\{8,9,10,11,12\}\) in subset \(B\) (shown in blue). Correspondingly, $5! \cdot 5! = 14400$ equivalent optimal sequences are detected.
Figure~\ref{fff1} reports the corresponding Gantt charts for Johnson’s sequence and for an
equivalent optimal sequence in which the orders of jobs \(\{3,4,5,6,7\}\) and
\(\{8,9,10,11,12\}\) are reversed.

\medskip

\noindent An even stronger result in terms of equivalent optimal schedules can be derived from
Properties~\ref{prop5} and~\ref{prop6} whenever either
\( P_{1} \leq P_{2} - p^{\max}_{B2} \) or
\( P_{2} \leq P_{1} - p^{\max}_{A1} \) holds.

\begin{property}
\label{prop5}
If \( P_{1} \leq P_{2} - p^{\max}_{B2} \), let \( k_A \) denote the smallest index such that
\[\sum_{i \in A: \, p_{1,i}\leq p_{1,k_A}} p_{1,i}
\leq
\sum_{i \in A: \, p_{1,i}\leq p_{1,k_A}} p_{2,i} - p^{\max}_{A1}
\]
and
\[
p_{1,k_A} < p_{1,k_A+1}.
\]
Then, an optimal schedule is obtained by first sequencing set $A$,
with jobs \(1,\dots,k_A\) in Johnson’s order and jobs \(k_A+1,\dots,n_A\) in any order,
followed by set $B$ sequenced in any order.
Moreover, if \( k_A \log k_A \leq n \), then an optimal solution to problem
$F2 \; || \; C_{\max}$ can be computed in linear time.
\end{property}

\begin{proof}
The sequencing of set $A$ follows directly from Property~\ref{prop1}.
For the second part, since
\[
P_{1} = P_{1A} + P_{1B} \leq P_{2} - p^{\max}_{B2}
= P_{2A} + P_{2B} - p^{\max}_{B2},
\]
we obtain
\[
P_{2A} - P_{1A} \geq P_{1B} - P_{2B} + p^{\max}_{B2}.
\]
Therefore, regardless of the sequence chosen for set $B$, the makespan is always given by
\( C_{2,n_A} + P_{2B} \), and no idle time is induced on machine $M_2$ while processing set
$B$.
The time-complexity argument follows exactly the same reasoning as in
Property~\ref{coro1}.
\end{proof}

\begin{property}
\label{prop6}
If \( P_{2} \leq P_{1} - p^{\max}_{A1} \), let \( k_B \) denote the largest index such that
\[\sum_{i \in B: \, p_{2,i}\leq p_{2,k_B}} p_{2,i}
\leq
\sum_{i \in B: \, p_{2,i}\leq p_{2,k_B}} p_{1,i} - p^{\max}_{B2}
\]
and
\[
p_{2,k_B} > p_{2,k_B-1}
\]
Then, an optimal schedule is obtained by first sequencing set $A$ in any order,
followed by jobs \(n_A+1,\dots,k_B-1\) in any order, and finally by jobs
\(k_B,\dots,n\) in Johnson’s order.
Moreover, if \( (n-k_B+1)\log(n-k_B+1) \leq n \), then an optimal solution to
problem $F2 \; || \; C_{\max}$ can be computed in linear time.
\end{property}

\begin{proof}
By the reversibility property (Property~\ref{prop11}), the conditions of
Property~\ref{prop6} correspond to the conditions of Property~\ref{prop5} applied to the
reversed instance.
\end{proof}

\medskip

\begin{remark}
In Properties~\ref{prop1} and~\ref{prop5}, there may exist several jobs
$h < k_A$ such that $p_{1,h} = p_{1,k_A}$.
Let \( k'_A \) denote the smallest index among these jobs.
Then, an optimal schedule can be obtained by sorting only jobs
\(1,\dots,k'_A-1\) and sequencing them in that order, followed by jobs
\(k'_A,\dots,k_A\) in any order (as they share the same processing time on machine $M_1$),
and finally by jobs \(k_A+1,\dots,n_A\) in any order.
A symmetric argument applies to Properties~\ref{prop2} and~\ref{prop6}:
there may exist several jobs $l > k_B$ such that $p_{2,l} = p_{2,k_B}$.
Let \( k'_B \) denote the largest index among these jobs.
Then, only jobs \(k'_B+1,\dots,n\) need to be sorted in subset $B$, preceded by jobs
\(k_B,\dots,k'_B\) and \(n_A+1,\dots,k_B-1\), both sequenced in any order.
\end{remark}

\medskip

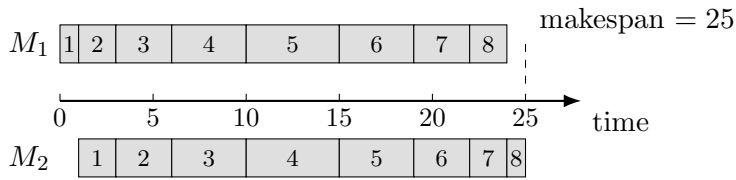
\begin{figure*}
	\begin{center}
		\begin{tabular}{cccccccccc}
			\toprule
			Job $j$   & 1 & 2 & 3 & 4 & 5 & 6 & 7 & 8 \\
			\midrule
			$p_{1,j}$ & 1 & 2 & 3 & 4 & 5 & 4 & 3 & 2 \\
			$p_{2,j}$ & 2 & 3 & 4 & 5 & 4 & 3 & 2 & 1 \\
			\bottomrule
		\end{tabular}
	\end{center}

	\medskip

	\begin{tikzpicture}[scale=0.7,x=0.35cm,y=0.9cm,>=Latex]
		\def\H{0.8}
		\def\Ygap{1.2}

		\draw[->,thick] (0,0) -- (28,0) node[below right] {time};
		\foreach \t in {0,5,10,15,20,25}
		\draw (\t,0) -- ++(0,0.15) node[below=2pt] {\small \t};

		\node[left] at (0, \Ygap) {\bf $M_1$};
		\node[left] at (0,-\Ygap) {\bf $M_2$};

		\draw[fill=gray!25,draw=black] (0,\Ygap-\H/2) rectangle ++(1,\H);
		\node at (0.5,\Ygap) {\footnotesize 1};

		\draw[fill=gray!25,draw=black] (1,\Ygap-\H/2) rectangle ++(2,\H);
		\node at (2,\Ygap) {\footnotesize 2};

		\draw[fill=gray!25,draw=black] (3,\Ygap-\H/2) rectangle ++(3,\H);
		\node at (4.5,\Ygap) {\footnotesize 3};

		\draw[fill=gray!25,draw=black] (6,\Ygap-\H/2) rectangle ++(4,\H);
		\node at (8,\Ygap) {\footnotesize 4};

		\draw[fill=gray!25,draw=black] (10,\Ygap-\H/2) rectangle ++(5,\H);
		\node at (12.5,\Ygap) {\footnotesize 5};

		\draw[fill=gray!25,draw=black] (15,\Ygap-\H/2) rectangle ++(4,\H);
		\node at (17,\Ygap) {\footnotesize 6};

		\draw[fill=gray!25,draw=black] (19,\Ygap-\H/2) rectangle ++(3,\H);
		\node at (20.5,\Ygap) {\footnotesize 7};

		\draw[fill=gray!25,draw=black] (22,\Ygap-\H/2) rectangle ++(2,\H);
		\node at (23,\Ygap) {\footnotesize 8};

		\draw[fill=gray!25,draw=black] (1,-\Ygap-\H/2) rectangle ++(2,\H);
		\node at (2,-\Ygap) {\footnotesize 1};

		\draw[fill=gray!25,draw=black] (3,-\Ygap-\H/2) rectangle ++(3,\H);
		\node at (4.5,-\Ygap) {\footnotesize 2};

		\draw[fill=gray!25,draw=black] (6,-\Ygap-\H/2) rectangle ++(4,\H);
		\node at (8,-\Ygap) {\footnotesize 3};

		\draw[fill=gray!25,draw=black] (10,-\Ygap-\H/2) rectangle ++(5,\H);
		\node at (12.5,-\Ygap) {\footnotesize 4};

		\draw[fill=gray!25,draw=black] (15,-\Ygap-\H/2) rectangle ++(4,\H);
		\node at (17,-\Ygap) {\footnotesize 5};

		\draw[fill=gray!25,draw=black] (19,-\Ygap-\H/2) rectangle ++(3,\H);
		\node at (20.5,-\Ygap) {\footnotesize 6};

		\draw[fill=gray!25,draw=black] (22,-\Ygap-\H/2) rectangle ++(2,\H);
		\node at (23,-\Ygap) {\footnotesize 7};

		\draw[fill=gray!25,draw=black] (24,-\Ygap-\H/2) rectangle ++(1,\H);
		\node at (24.5,-\Ygap) {\footnotesize 8};

		\draw[dashed] (25,1.2) -- (25,0.2);
		\node[above] at (31,1.2)  { makespan $=25$};
	\end{tikzpicture}
	\caption{An instance with unique optimal sequence}.
	\label{f2}
\end{figure*}

\noindent
As a counterpart, we show that in the worst case a complete sorting of both subsets
$A$ and $B$ may be necessary.
To this end, consider a two-machine flow shop instance with $2n$ jobs.
The first $n$ jobs belong to subset $A$ and have processing times
$p_{1,j}=j$ and $p_{2,j}=j+1$ for all $j=1,\dots,n$.
The remaining $n$ jobs belong to subset $B$ and have processing times
$p_{1,j}=2n-j+2$ and $p_{2,j}=2n-j+1$ for all $j=n+1,\dots,2n$.
It is straightforward to verify that the unique optimal schedule is given by the
sequence
\[
1,2,\dots,n,n+1,\dots,2n,
\]
and therefore a complete sorting of subset $A$ in nondecreasing order of processing
times on machine $M_1$ and of subset $B$ in nonincreasing order of processing times on
machine $M_2$ is required.
Figure~\ref{f2} illustrates this instance and the corresponding unique optimal schedule
for $n=4$, with makespan equal to $25$.
\bigskip

\section{Probabilistic analysis}
\label{prob}
In this section, $p^{\max}$ denotes the upper bound of the uniform distribution (i.e., the maximum possible processing time), consistently with the experimental setting.
Here, we present a probabilistic analysis for instances with $n$ jobs, where
processing times are independently and uniformly distributed over the set
$\{1,\dots,p^{\max}\}$.
We focus on parameter $k_A$ and subset $A$.
By the reversibility property (Property~\ref{prop11}), the same analysis applies
symmetrically to parameter $k_B$ and subset $B$.

Our goal is to estimate the probability that the value of $k_A$ is sufficiently small
to guarantee linear-time complexity, and to analyze how this probability behaves as $n$
grows.
Unless otherwise stated, we assume that $p^{\max}=\Theta(n)$, which is the standard
setting in flow shop benchmark instances.

\subsection{Step 1: Lower bound on the size of set $A$}

We fix a parameter $\alpha$ with $0 < \alpha < 1$.
For each job $i$, we consider the following restrictive condition:
\[
C_1(i) :=
\bigl(p_{1,i} \le \alpha p^{\max} \bigr)
\;\land\;
\bigl(p_{2,i} > \alpha p^{\max} \bigr).
\]
Clearly, condition $C_1(i)$ implies that job $i$ belongs to set $A$, although the converse
does not necessarily hold.

\begin{lemma}[Probability of condition $C_1$]
For any job $i$, the probability that condition $C_1(i)$ holds is
\[
\Pr[C_1(i)] = \alpha (1-\alpha).
\]
\end{lemma}

\begin{proof}
Since $p_{1,i}$ and $p_{2,i}$ are independent and uniformly distributed over
$\{1,\dots,p^{\max}\}$, we have
\[
\Pr(p_{1,i} \le \alpha p^{\max}) = \alpha, 
%\quad\text{and}\quad
\Pr(p_{2,i} > \alpha p^{\max}) = 1-\alpha.
\]
Correspondingly, $\Pr[C_1(i)] = \alpha (1-\alpha)$.
\end{proof}

Let $X$ be the random variable counting the number of jobs satisfying condition $C_1$.
Then $X$ follows a binomial distribution with parameters $n$ and
$p=\alpha(1-\alpha)$.
By the binomial formula,
\[
\Pr(X = k) = \binom{n}{k} [\alpha(1-\alpha)]^k [1-\alpha(1-\alpha)]^{n-k}.
\]

Since we are interested in the probability that at least $k_A$ jobs satisfy condition $C_1$, we
define
\[
P_1^* :=
\Pr(X \ge k_A)
= \]
\[ =
\sum_{k=k_A}^{n}
\binom{n}{k} [\alpha(1-\alpha)]^k [1-\alpha(1-\alpha)]^{n-k}.
\]

\subsection{Step 2: Lower bound on the sum of processing time differences}

From Step~1, it follows that with probability lower bounded by $P_1^*$ there exist at least
$k_A$ jobs $i \in A$ such that $p_{1,i} \le \alpha p^{\max}$.
Consider the first $k_A$ jobs in set $A$ according to Johnson’s sequence.
The inequality $p_{1,i} \le \alpha p^{\max}$ necessarily holds for all these jobs.
For each such job $i$, let $\delta_i := p_{2,i} - p_{1,i}$.

Conditioned on a fixed value of $p_{1,i}$, the random variable $\delta_i$ is uniformly
distributed over
\[
\{0,\dots,p^{\max}-p_{1,i}\}.
\]

Let $p' := (1-\alpha)p^{\max}$.
We further restrict the support of each $\delta_i$ to the smaller set
$\{0,\dots,p'\}$.
%This restriction yields a valid lower bound, since we condition on an event of positive probability and then further restrict the support.

This restriction yields a valid lower bound because we replace the original event with a more restrictive one:
we require each $\delta_i$ to lie in $\{0,\dots,p'\}$, and we then compute the probability of
$\sum_{i=1}^{k_A}\delta_i \geq p^{\max}$ under this restriction.
Correspondingly,
the variables $\delta_1,\dots,\delta_{k_A}$ are independent and
uniformly distributed over $\{0,\dots,p'\}$.

We can then compute a lower bound on the probability that
\[
\sum_{i=1}^{k_A} \delta_i \ge p^{\max}
\]
by means of dynamic programming.

\medskip

Let $P(j,\gamma)$ denote the probability that $\sum_{i=1}^{j} \delta_i = \gamma$, i.e., the ratio between the number of realizations yielding exactly $\gamma$ and $(p'+1)^j$, where $(p'+1)^j$ represents the total number of possible realizations.

The maximum attainable value of $\gamma$ is $j p'$.
The dynamic program is defined as follows:
\begin{itemize}
\item Initial conditions: \\
$
P(0,0) = 1,
\;
P(0,\gamma) = 0 \; \forall \gamma > 0$.
\item Boundary conditions: \\
$P(j,\gamma) = 0 \quad \text{if } \gamma < 0$.

\item Recursion: \\
%P(j,\gamma)
%=
%\frac{1}{p'+1}
%\sum_{i=0}^{p'} P(j-1,\gamma-i),
%j=1,\dots,k_A.
$P(j,\gamma)
=
\frac{\sum_{i=0}^{p'} P(j-1,\gamma-i)}{p'+1}, \;
j=1,\dots,k_A$.
\end{itemize}

Finally, the probability that $\sum_{i=1}^{k_A} \delta_i \ge p^{\max}$
holds is given by
\[
P_2^*
=
1 - \sum_{\gamma=0}^{p^{\max}-1} P(k_A,\gamma).
\]

Correspondingly, the product $P^* = P_1^* \cdot P_2^*$
provides a lower bound on the probability that just the first $k_A$ jobs need to be
sorted in subset $A$.
By the reversibility property (Property~\ref{prop11}), an analogous bound holds for subset
$B$ and parameter $k_B$.

In Table~\ref{tabprob1}, we consider the standard benchmark setting $p^{\max}=n$, and for each value of $n$  we select the pair $(k_A,\alpha)$
that maximizes the bound
$P^* = P_1^* \cdot P_2^*$
%\[
%P^*(n,p^{\max};k_A,\alpha)=P_1^*(n;k_A,\alpha)\cdot P_2^*(p^{\max};k_A,\alpha),
%\]
subject to the complexity constraint $k_A \log k_A \le p^{\max}$.
The parameter $\alpha$ is searched over the grid $\{0.01,0.02,\dots,0.99\}$.
We report the resulting $k_A$, $\alpha$, and the corresponding value of $P^*$.
The table indicates that $P^* > 1 - 10^{-1}$ for $n=40$, $P^* > 1 - 10^{-3}$ for $n=100$,
and $P^* > 1 - 10^{-9}$ for $n=200$, supporting the claim that the proposed approach runs in overall linear time with high probability
under uniformly distributed processing times.

\begin{table}[ht!]
\begin{center}
{\footnotesize
\begin{tabular}{rrrr}
			$n$   & $k_A$ & $\alpha$  & $P^*$  \\ \hline
$20$ & $4$ & $0.35$ & $ \approx 0.547$ \\ 
$40$ & $6$ & $0.40$ & $\approx 0.912$ \\ 
$60$ & $7$ & $0.38$ & $\approx 0.987$ \\ 
$80$ & $9$ & $0.42$ & $\approx 0.998$ \\ 
$100$ & $10$ & $0.41$ & $> 1 - 10^{-3}$ \\ 
$120$ & $12$ & $0.45$ & $> 1 - 10^{-4}$ \\ 
$140$ & $13$ & $0.45$ & $> 1 - 10^{-5}$ \\ 
$160$ & $14$ & $0.43$ & $> 1 - 10^{-6}$ \\ 
$180$ &  $15$ & $0.42$ &  $> 1 - 10^{-7}$ \\  
$200$ &  $17$ &  $0.47$ &   $> 1 - 10^{-9}$ \\ \hline
\end{tabular}
}
\end{center}
\caption{Values of $P^*$ for increasing value of $n$}
\label{tabprob1}
\end{table}

%\begin{remark}
%For fixed $k_A$ and $\alpha$,  $P_1^*$ is non-decreasing in $n$,
%since increasing $n$ corresponds to adding independent trials.
%\end{remark}

\section{Computational testing}
\label{compres}

In this section, we present a computational analysis of the sizes of the job subsets
that need to be sorted.
Specifically, we focus on the values $k_A$ and $k'_A$ (with $k'_A \le k_A$) for subset $A$,
and on the values
$\overline{k_B}=n-k_B+1$ 
(note that $\overline{k_B} = | \{k_B, \dots, n\}|$)
and $\overline{k'_B}=n-k'_B+1$ (with
$\overline{k'_B}\le \overline{k_B}$) for subset $B$.
In standard flow shop benchmark instances, processing times are typically generated
according to uniform distributions; see, e.g., \cite{Taillard}.
Accordingly, we considered instances of size
$S \in \{10^2,10^3,10^4,10^5\}$, where processing times were independently generated
according to uniform distributions on $[1,S]$ and $[1,10S]$.
For each pair (size, distribution), $100$ instances were generated, for a total of
$800$ instances.
All algorithms were implemented in \texttt{C++20}, and the experiments were run on a
PC equipped with an AMD Ryzen~7~5800H CPU clocked at 4.4~GHz.

The results are summarized in Table~\ref{tabresUniform}.
The first column reports the number of jobs.
The second column indicates the maximum processing time $p^{\max}$ of the corresponding
uniform distribution.
The third and fourth columns report the maximum values of $k'_A$ and $k_A$,
respectively, observed among the instances with same size and same distribution.
The fifth and sixth columns report the maximum values of $\overline{k'_B}$ and
$\overline{k_B}$, respectively.
The seventh column reports, for each (size, distribution), the number of instances in
which either Property~\ref{prop5} or Property~\ref{prop6} holds.
The eighth column reports the average computational time (in seconds) of our algorithm.
Finally, the last column reports the ratio between the average CPU time required by a standard
implementation of Johnson’s rule with full sorting of sets $A$ and $B$ and the average CPU time
required by our approach.

\begin{table*}[ht!]
	\centering
	\setlength{\tabcolsep}{3.5pt}
	% \fontsize{10}{8}\selectfont
{\footnotesize
	\begin{tabular}{ccccccccc}
		\toprule
		n      & $p^{\max}$ & $k'_A$ & $k_A$ & $\overline{k'_B}$ & $\overline{k_B}$ & Prop5,6 & $t_{avg}$          & $\tau$ \\
\midrule
100 & 100 & 6 & 3 & 7 & 5 & 83 & $1.8\times10^{-6}$ & 2.27 \\
100 & 1000 & 5 & 3 & 7 & 6 & 87 & $1.6\times10^{-6}$ & 2.45 \\
1000 & 1000 & 6 & 4 & 6 & 4 & 94 & $2.8\times10^{-6}$ & 15.05 \\
1000 & 10000 & 6 & 4 & 5 & 4 & 95 & $2.8\times10^{-6}$ & 14.64 \\
10000 & 10000 & 6 & 5 & 6 & 4 & 97 & $1.1\times10^{-5}$ & 45.36 \\
10000 & 100000 & 6 & 5 & 4 & 3 & 97 & $1.1\times10^{-5}$ & 45.28 \\
100000 & 100000 & 8 & 6 & 6 & 3 & 100 & $9.9\times10^{-5}$ & 63.79 \\
100000 & 1000000 & 7 & 6 & 6 & 5 & 100 & $1.0\times10^{-4}$ & 64.40 \\
\bottomrule

	\end{tabular}
}	\caption{Testing instances with uniform distribution}
	\label{tabresUniform}
\end{table*}

\medskip

\noindent
Table~\ref{tabresUniform} shows that, in all tested instances, both $k_A$ and $k_B$ never
exceeded $8$, meaning that the sorting required for both subsets $A$ and $B$ was limited
to at most $8$ jobs.
Consequently, for all instances, an optimal solution of the corresponding
$F2 \; || \; C_{\max}$ problem could be computed in linear time.
Moreover, for the majority of the instances ($753$ out of $800$), either
Property~\ref{prop5} or Property~\ref{prop6} holds, implying that only one of the two
subsets required any sorting.

%\pagebreak

In addition to the uniform distribution on $[1,p^{\max}]$, we also considered three other
discrete probability distributions commonly used to model variability in processing
times.
For each distribution, the parameters were chosen so that the expected processing time
was approximately $p^{\max}/2$, ensuring comparability across all distributions.

\begin{itemize}
	\item \textbf{Negative binomial distribution.}
	      We use the $(r,p)$ parameterization, where $r>0$ denotes the number of successes
	      and $0<p<1$ the probability of success.
	      The mean is $\mathbb{E}[X]=r(1-p)/p$.
	      We fix $r=5$ and set $p = \frac{r}{r+\mu}$ with $\mu = p^{\max}/2$, yielding a
	      moderately dispersed distribution centered around $p^{\max}/2$.

	\item \textbf{Geometric distribution.}
	      A geometric distribution with parameter $p$ has mean
	      $\mathbb{E}[X]=(1-p)/p$.
	      Setting $p = \frac{2}{p^{\max}+2}$ yields
	      $\mathbb{E}[X]\approx p^{\max}/2$.
	      This distribution is highly skewed and emphasizes short processing times.

	\item \textbf{Poisson distribution.}
	      A Poisson distribution with rate $\lambda$ has mean
	      $\mathbb{E}[X]=\lambda$.
	      We set $\lambda = p^{\max}/2$, yielding a discrete distribution with lighter
	      tails than the geometric and negative binomial distributions.
\end{itemize}
%\pagebreak

\begin{figure}[ht]
	\centering
	\begin{minipage}{.4\textwidth}
		\centering
		\includegraphics[width=\textwidth]{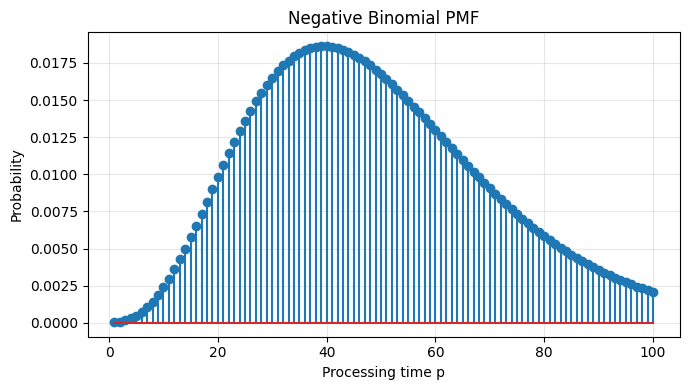}
		\captionof{figure}{ \footnotesize{Negative Binomial distribution PMF}}
		\label{fig:PMF_neg_bin}
	\end{minipage}%
	\hspace{2mm}
	\begin{minipage}{.4\textwidth}
		\centering
		\includegraphics[width=\textwidth]{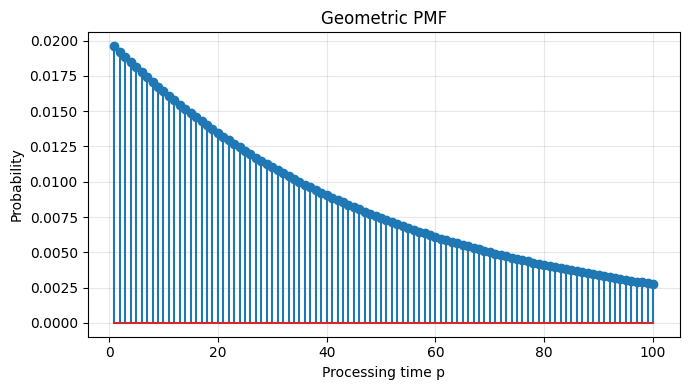}
		\captionof{figure}{\footnotesize{Geometric distribution PMF}}
		\label{fig:PMF_geo}
	\end{minipage}
	\hspace{2mm}
	\begin{minipage}{.4\textwidth}
		\centering
		\includegraphics[width=\textwidth]{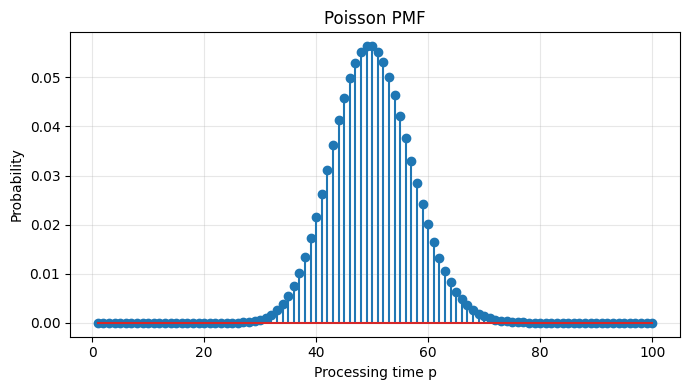}
		\captionof{figure}{\footnotesize{Poisson distribution PMF}}
		\label{fig:PMF_poi}
	\end{minipage}
\end{figure}

Figures~\ref{fig:PMF_neg_bin}, \ref{fig:PMF_geo}, and~\ref{fig:PMF_poi} illustrate the probability mass functions (PMF) of the considered distributions.

The results are summarized in Table~\ref{tabresAllDistrib}, which shows that the choice
of distribution has a negligible impact on the values of $k_A$ and $k_B$.
These experiments were conducted with $n=p^{\max}=1000$.

\begin{table*}[ht]
	\centering
	\setlength{\tabcolsep}{3.5pt}
{\footnotesize
	\begin{tabular}{lrrrrr}
		\toprule
		Distribution      & $k'_A$ & $k_A$ & $\overline{k'_B}$ & $\overline{k_B}$ & Prop5,6 \\
		\midrule
		geometric            & 4        & 6      & 3                    & 5                 & 96     \\
		negative binomial & 5        & 6      & 4                    & 7                 & 96      \\
		poisson               & 3        & 7      & 5                    & 9                 & 96       \\
		\bottomrule
	\end{tabular}
}	\caption{Testing instances with different distributions where $n=p^{\max}=1000$}
	\label{tabresAllDistrib}
\end{table*}

%\pagebreak

\section{Conclusions}
\label{conclus}

We revisited Johnson’s rule for makespan minimization in the two-machine flow shop
scheduling problem.
We provided a sufficient condition, verifiable in linear time, to determine whether
an optimal solution can be computed in linear time.
We also showed that, although instances exist for which a complete sorting of the
subsets $A$ and $B$ is required, in practice, under the standard uniform distribution,
only a very limited amount of sorting is needed.
As a consequence, an optimal sequence can be computed in linear time with high
probability (lower bounded by $1-10^{-9}$ for instances of size $n = 200$).

We also remark that all permutations of the jobs in $A$ and $B$ that do not need to be
sorted can be part of an optimal schedule, thus inducing correspondingly large sets
of equivalent optimal solutions.

Although it is beyond the scope of this paper, a tailored application of the proposed
properties appears suitable for the rapid detection of subsets of distinct optimal
sequences, in the spirit of \cite{Briand2006} and \cite{Hadda2018}.
For instance, the application of Property~\ref{prop5} to the example proposed
in \cite{Briand2006} immediately detects $48$ distinct optimal sequences, which is
significantly more than the $13$ sequences reported in \cite{Briand2006} and matches
the value reported in \cite{Hadda2018}.
Moreover, when a straightforward partial enumeration is combined with a standard
application of the proposed properties, $492$ distinct optimal sequences are detected, a
value that is larger than the $324$ optimal sequences reported in \cite{Hadda2018}.
Investigating these aspects in greater depth constitutes an interesting direction
for future research.

%\noindent
%{\bf Data availability}

%Data will be made available on request.

%\footnotesize

\end{document}